\documentclass[12pt]{amsart}
\usepackage{mathrsfs}
\usepackage{amssymb}

\usepackage{titletoc}
\input xy
  \xyoption{all}
\usepackage{amscd}
\usepackage{amsmath, amssymb}
\usepackage{amsfonts}
\usepackage[colorlinks,linkcolor=blue,citecolor=blue, pdfstartview=FitH]{hyperref}

\numberwithin{equation}{section}

\theoremstyle{plain}
\newtheorem{thm}{Theorem}[section]
\newtheorem{theorem}[thm]{Theorem}
\newtheorem{lemma}[thm]{Lemma}
\newtheorem{corollary}[thm]{Corollary}
\newtheorem{proposition}[thm]{Proposition}

\theoremstyle{definition}

\newtheorem{question}[thm]{Question}

\newtheorem{remark}[thm]{Remark}

\newtheorem{definition}[thm]{Definition}

\newtheorem{example}[thm]{Example}

\newtheorem{defn-thm}[thm]{Definition-Theorem}

\newcommand{\C}{{\mathbb C}}

\renewcommand{\S}{{\mathbb S}}

\newcommand{\qtq}[1]{\quad\mbox{#1}\quad}
\newcommand{\bp}{\bar{\partial}}
\newcommand{\Om}{\Omega}

\newcommand{\ts}{\otimes}

\newcommand{\btheorem}{\begin{theorem}}
\newcommand{\etheorem}{\end{theorem}}
\newcommand{\bproposition}{\begin{proposition}}
\newcommand{\eproposition}{\end{proposition}}
\newcommand{\bdefinition}{\begin{definition}}
\newcommand{\edefinition}{\end{definition}}
\newcommand{\bcorollary}{\begin{corollary}}
\newcommand{\ecorollary}{\end{corollary}}
\newcommand{\bproof}{\begin{proof}}
\newcommand{\eproof}{\end{proof}}
\newcommand{\bremark}{\begin{remark}}
\newcommand{\eremark}{\end{remark}}
\newcommand{\eexample}{\end{example}}
\newcommand{\bexample}{\begin{example}}
\newcommand{\la}{\langle}
\newcommand{\elemma}{\end{lemma}}
\newcommand{\blemma}{\begin{lemma}}
\newcommand{\ra}{\rangle}
\newcommand{\sq}{\sqrt{-1}}

\newcommand{\p}{\partial}

\renewcommand{\bar}{\overline}

\renewcommand{\phi}{\varphi}

\newcommand{\ee}{\end{eqnarray*}}
\newcommand{\be}{\begin{eqnarray*}}

\newcommand{\beq}{\begin{equation}}
\newcommand{\eeq}{\end{equation}}

\newcommand{\bd}{\begin{enumerate}}
\newcommand{\ed}{\end{enumerate}}

\renewcommand{\>}{\rightarrow}

\newcommand{\bbe}{\bar{\beta}}

\newcommand{\al}{{\alpha}}
\newcommand{\ga}{{\gamma}}

\usepackage{fancyhdr}
\renewcommand{\leftmark}{{\tiny RC-positivity and rigidity of harmonic maps into Riemannian manifolds}}

\begin{document}
\title{RC-positivity and rigidity of harmonic maps into Riemannian manifolds} \makeatletter
\let\uppercasenonmath\@gobble
\let\MakeUppercase\relax
\let\scshape\relax
\makeatother

\author{Jun Wang}
\date{}
\address{{ Address of Jun Wang:  Tsinghua University, Yau Mathematical Science Center, Beijing, 100084, China.}}

\author{Xiaokui Yang}
\date{}
\address{{Address of Xiaokui Yang: Morningside Center of Mathematics, Institute of
        Mathematics, Hua Loo-Keng Center of Mathematical Sciences,
        Academy of Mathematics and Systems Science,
        Chinese Academy of Sciences, Beijing, 100190, China.}}
\email{\href{mailto:xkyang@amss.ac.cn}{{xkyang@amss.ac.cn}}}

\noindent\thanks{This work was partially supported   by China's
Recruitment
 Program of Global Experts
 }
\maketitle \vskip -0.5cm {\hfill\small{\emph{Dedicated to Professor
Lo Yang on the occasion of his 80th birthday \hfill}}

\begin{abstract} In this paper, we show that every harmonic map from   a compact K\"ahler manifold  with
uniformly RC-positive curvature  to a Riemannian manifold with
non-positive complex sectional curvature  is constant. In
particular, there is no non-constant harmonic map from a compact
K\"ahler manifold with positive holomorphic sectional curvature to a
Riemannian manifold with non-positive complex sectional curvature.

\end{abstract}
{{\setcounter{tocdepth}{1} \tableofcontents}}

\section{Introduction}

The classical result of Eells-Sampson \cite{ES64} asserts that
harmonic maps from  compact Riemannian manifolds with positive Ricci
curvature to Riemannian manifolds with non-positive sectional
curvature must be constant. By using Schwarz calculation and maximum
principle, S.-T. Yau established in \cite{Yau78} that holomorphic
maps from compact K\"ahler manifolds with positive Ricci curvature
to Hermitian manifolds with non-positive holomorphic bisectional
curvature must be constant. Recently, the second author extended
Yau's result and obtained in \cite{Yang18a} that  holomorphic maps
from compact Hermitian manifolds with
 positive holomorphic sectional curvature to Hermitian manifolds with
non-positive holomorphic bisectional curvature must be constant, and
one of the key ingredients in the proof is a notion called
``RC-positivity" introduced in \cite{Yang18}. Moreover, Y.-T. Siu
achieved in \cite{Siu80} the complex analyticity of harmonic maps
into compact K\"ahler manifolds with strongly negative curvature.
Later, Jost-Yau explored in \cite{JY91, JY93} rigidity theorems  for
various harmonic maps from Hermitian manifolds into Riemannian
manifolds with non-positive complex sectional curvature (see also
\cite{Sam85, Sam86}).
There are many generalizations along this line, and we refer to \cite{EL78, Siu82, EL83, Ohn87, EL88, Uda88, OU90, YZ91, MSY93, Uda94, Ni99, Tos07, Zhang12, YHD13, Dong13, LY14, YZ18, Yang18a, Yang18b, Ni18} and the references therein. \\

 In this paper, we aim to extend the above rigidity results to harmonic maps
 (and pluri-harmonic maps) from compact complex manifolds to
 Riemannian manifolds by using RC-positivity and
 ideas in \cite{ES64, Yau78, Siu80, Sam85, Sam86, JY91, JY93, LY14, Yang18, Yang18a,
 Yang18b}. More precisely, we consider harmonic maps (resp. pluri-harmonic maps) from a compact complex manifold $M$ with uniformly RC-positive curvature to a Riemannian manifold $N$
 with non-positive complex sectional curvature. Let's recall these
 curvature notions. In \cite{Yang18b}, the second author introduced  the following
definition.

\bdefinition
 A Hermitian holomorphic vector bundle $(\mathscr E,h)$ over a
complex manifold $M$ is called \emph{uniformly RC-positive}  if at
every point $q\in M$, there exists some vector $u\in T^{1,0}_qM$
such that for any nonzero vector $v\in \mathscr E_q$, one has  \beq
R^\mathscr E(u,\bar u,v,\bar v)>0.\eeq \edefinition

\noindent There are many complex manifolds which have uniformly
RC-positive tangent bundles (i.e., the tangent bundle admits a
Hermitian metric such that its Chern curvature is uniformly
RC-positive). For instances,

 \bd  \item  K\"ahler manifolds with positive holomorphic sectional curvature
(\cite[Theorem~5.1]{Yang18b});

\item  Hopf manifolds $\S^1\times \S^{2n+1}$ with standard metrics
(\cite[formula~(6.4)]{LY14});

\item products of complex manifolds with uniformly RC-positive
tangent bundles. \ed For more details, we refer to \cite{Yang18,
Yang18a, Yang18b} and the references therein.\\

 As analogous to Siu's strong negativity \cite{Siu80},
Sampson
 proposed  in \cite{Sam85} the following definition (see also
\cite[Definition~4.2]{JY91}):

\bdefinition \label{samdef}Let $(N,g)$ be a  Riemannian manifold.
The (complexified) curvature tensor $R$ of $(N,g)$ is said to have
\emph{non-positive complex sectional curvature} if \beq R(Z,\bar W,
W,\bar Z)\leq 0 \label{cscdef}\eeq for any $Z,W\in T_\C
N$.\edefinition

\noindent If $(N,g)$ has non-positive complex sectional curvature,
then it has non-positive Riemannian sectional curvature. Moreover,
 a K\"ahler manifold $(N,g)$ has
strongly non-positive curvature in the sense of Siu if and only if
its background Riemannian metric has non-positive complex sectional
curvature (e.g. \cite[Theorem~4.4]{LSYY17}). Hence, compact
quotients of bounded symmetric domains of type $\mathrm{I}_{mn} (mn
\geq2)$, $\mathrm{II}_n (n\geq  3)$, $\mathrm{III}_n (n\geq 2)$, and
$\mathrm{IV}_n (n \geq 3)$ all have non-positive complex sectional
curvature.

\noindent  The first main result of this paper is the following.

 \btheorem
\label{main0} Let $(M,h)$ be a compact K\"ahler manifold with
uniformly RC-positive curvature and $(N,g)$ be a Riemannian manifold
with non-positive complex sectional curvature.  Then  there is no
non-constant harmonic map from $(M,h)$ to $(N,g)$. \etheorem

\noindent By Theorem \ref{main0} and \cite[Theorem~5.1]{Yang18b}, we
obtain

\bcorollary\label{coro} There is no non-constant harmonic map from a
compact K\"ahler manifold with positive holomorphic sectional
curvature to a Riemannian manifold with non-positive complex
sectional curvature.

\ecorollary

\bremark We proved in \cite[Corollary~1.3]{Yang18b} that for a
compact K\"ahler manifold $M$, if it  admits a Hermitian metric with
uniformly RC-positive curvature, then $M$ is projective and
rationally connected. By using algebraic geometry,  we can prove a
more general result that: \emph{every harmonic map from a rationally
connected manifold to a Riemannian manifold with non-positive
complex sectional curvature is constant.} It is well-known that Fano
manifolds are all rationally connected (\cite{Cam92, KMM92}).
 Hence,
Theorem \ref{main0} and Corollary \ref{coro} are  generalizations of
the following well-known result: every harmonic map from a Fano
manifold (i.e. a compact K\"ahler manifold with positive Ricci
curvature \cite{Yau78a}) to a Riemannian manifold with non-positive
(complex) sectional curvature is constant. \eremark

In the following, we give the proof of Theorem \ref{main0} in a more
general setting.  Let $f:(M,h)\>(N,g)$ be a smooth map from a
compact Hermitian manifold $(M,h)$ to a Riemannian manifold $(N,g)$.
Jost and Yau proposed in \cite[(H2)]{JY93} the following concept:
$f$ is called a \emph{Hermitian harmonic map} if it satisfies \beq
h^{\alpha\bar\beta}\left(\frac{\p^2 f^i}{\p z^\alpha\p\bar
z^\beta}+\Gamma_{jk}^i\frac{\p f^j}{\p \bar z^\beta}\frac{\p f^k}{\p
z^\alpha}\right)\ts e_i=0. \eeq  It is easy to see that when $M$ is
K\"ahler, Hermitian harmonic maps are usual harmonic maps. On the
other hand, Jost-Yau proved in \cite{JY93} that  there always exist
Hermitian harmonic maps
 into Riemannian manifolds with
 non-positive sectional curvature, which generalizes the classical
 result of Eells-Sampson \cite{ES64}. Moreover, $f$ is called \emph{pluri-harmonic} if \beq
\left(\frac{\p^2 f^i}{\p z^\alpha\p\bar
z^\beta}+\Gamma_{jk}^i\frac{\p f^j}{\p \bar z^\beta}\frac{\p f^k}{\p
z^\alpha}\right)dz^\alpha\wedge d\bar z^\beta\ts
e_i=0.\label{0pluri}\eeq  Pluri-harmonic maps  are generalizations
of holomorphic maps.  Indeed,  a holomorphic map from a complex
manifold $M$ to a K\"ahler manifold $(N,g)$ is pluri-harmonic. The
second main result  is on rigidity of pluri-harmonic maps, which is
analogous to  rigidity of holomorphic maps obtained in
\cite{Yang18a,Yang18b}.

 \btheorem\label{main1} Let
$f:M\>(N,g)$ be a pluri-harmonic map from a  compact complex
manifold $M$ to a Riemannian manifold $(N,g)$ with non-positive
complex sectional curvature. If $T^{1,0}M$ is uniformly RC-positive,
then $f$ is a constant map. \etheorem

\noindent The key ingredient in the proof of Theorem \ref{main1} is
the following Chern-Lu type identity for pluri-harmonic maps.
\bproposition Let $f:(M,h)\>(N,g)$ be a pluri-harmonic map from a
 Hermitian manifold $(M,h)$ to a Riemannian manifold $(N,g)$.
Then we have the identity $$ \p\bp u =\la\p_{\mathscr F} \p f,
\p_{\mathscr F}\p f \ra + \left( R^M_{\alpha\bar \beta \gamma\bar
\delta} h^{\gamma \bar \nu} h^{\mu\bar \delta} g_{ij} f^{i }_\mu
{f^{j}_{\bar\nu}} - R^N_{i\ell k j}f^i_{\al}f^j_{\bbe}
\left(h^{\gamma\bar\delta}f^k_{\ga}f^\ell_{\bar\delta}\right)\right)dz^\alpha\wedge
d\bar z^\beta,$$
 where $u=|\p f|^2$, $R^M_{\alpha\bar \beta \gamma\bar
\delta}$ is the curvature of the Chern connection on $T^{1,0}M$ and
$\p_{\mathscr F}$ is the induced $(1,0)$ connection on the vector
bundle $\mathscr F=T^{*1,0}M\ts f^*T_N$. \eproposition

\noindent
  By using Theorem \ref{main1} and ideas in \cite{Siu80, Sam85, Sam86, YZ91, JY91, JY93,
 LY14}, we obtain the following result.  Note that, Theorem \ref{main0} is a
special case of it.

\btheorem\label{main2} Let $(M,h)$ be a compact  astheno-K\"ahler
manifold (i.e. $\p\bp\omega^{m-2}_h=0$) and $(N,g)$ be a  Riemannian
manifold with
 non-positive complex sectional  curvature. Let $f:(M,h)\>(N,g)$ be a Hermitian harmonic map. If $M$ admits a Hermitian metric $h_0$  with uniformly RC-positive curvature, then $f$ is a constant map.
 \etheorem

\noindent There are many astheno-K\"ahler manifolds which are
non-K\"ahler. For instances, non-K\"ahler complex surfaces and
Calabi-Eckmann manifolds $\S^{2p+1}\times \S^{2q+1}$
(\cite{Mat09}).\\

 As motivated by Theorem \ref{main0}, we propose:

\begin{question} Does there exist a non-constant Hermitian harmonic map (resp. harmonic map) from   a compact
Hermitian manifold  with (uniformly) RC-positive curvature  to a
Riemannian manifold with non-positive (complex) sectional curvature?

\end{question}

\noindent As a special case of Theorem \ref{main2}, one gets

\bcorollary Let $M$ be a compact complex surface and $(N,g)$ be a
Riemannian manifold with
 non-positive complex sectional  curvature. Let $f:(M,h)\>(N,g)$ be a Hermitian harmonic map.
 If $M$ admits a Hermitian metric $h_0$ with uniformly RC-positive curvature, then $f$ is a constant map.
 \ecorollary

\noindent As we pointed out before, the standard  Hopf surface
$\S^1\times \S^3$ admits uniformly RC-positive metrics, and so we
get

\bcorollary  There is no non-constant Hermitian harmonic map from
Hopf surface $\S^1\times \S^3$ (with an arbitrary Hermitian metric)
to Riemannian manifolds with
 non-positive complex sectional  curvature.
\ecorollary

\noindent \textbf{Acknowledgements.} The first author would like to
thank his advisor Professor Jian Zhou for his guidance. The second
author is very grateful to Professor
 K.-F. Liu and Professor S.-T. Yau for their support, encouragement and stimulating
discussions over  years. We would also like to thank Y.-X. Li and V.
Tosatti for some helpful discussions.

\vskip 2\baselineskip

\section{Harmonic maps from Hermitian manifolds to Riemannian manifolds}

\subsection{Connections on vector bundles}
Let $\mathscr E$ be a Hermitian {complex} vector bundle or a
Riemannian real vector bundle over a compact Hermitian manifold
$(M,h)$ and $\nabla^\mathscr E$ be a metric connection on $\mathscr
E$. There is a natural decomposition $\nabla^\mathscr
E=\nabla^{'\mathscr E}+\nabla^{''\mathscr E}$ where \beq
\nabla^{'\mathscr E}:\Gamma(M,\mathscr E)\>\Om^{1,0}(M,\mathscr E)
\qtq{and} \nabla^{''\mathscr E}:\Gamma(M,\mathscr
E)\>\Om^{0,1}(M,\mathscr E).
  \eeq
Moreover, $\nabla^{'\mathscr E}$ and $\nabla^{''\mathscr E}$  induce
two differential operators. The first one is
 $
\p_\mathscr E:\Om^{p,q}(M,\mathscr E)\>\Om^{p+1,q}(M,\mathscr E) $
which is defined by \beq \p_\mathscr E(\phi\ts
s)=\left(\p\phi\right)\ts s+(-1)^{p+q}\phi\wedge \nabla^{'\mathscr
E}s \eeq for any $\phi\in \Om^{p,q}(M)$ and $s\in \Gamma(M,\mathscr
E)$. The operator $ \bp_\mathscr E:\Om^{p,q}(M,\mathscr
E)\>\Om^{p,q+1}(M,\mathscr E) $ is defined similarly.   The operator
$\p_\mathscr E\bp_\mathscr E+\bp_\mathscr E\p_\mathscr E$ is
represented by the curvature tensor $R^\mathscr E \in \Gamma(M,
\Lambda^{1,1}T^*M\ts \mathscr E^*\ts \mathscr E)$.  The adjoint
operators of $\p,\bp,\p_\mathscr E$ and $\bp_\mathscr E$ are denoted
by $\p^*,\bp^*, \p^*_\mathscr E$ and $\bp_\mathscr E^*$
respectively.

\subsection{Harmonic maps between Riemannian manifolds}

 Let $(M,h)$ and $(N,g)$ be two compact Riemannian
 manifolds. Let $f:(M,h)\>(N,g)$ be a smooth map. If $\mathscr E=f^*(TN)$,
 then  $df$ can be regarded as  an $\mathscr E$-valued one form. There
 is an induced connection $\nabla^\mathscr E$ on $\mathscr E$ by the Levi-Civita connection on
 $TN$. In the local coordinates $\{x^\alpha\}_{\alpha=1}^m$, $\{y^i\}_{i=1}^n$
 on
 $M$ and $N$ respectively, the local frames of $\mathscr E$ are denoted by
 $e_i=f^*\left(\frac{\p}{\p y^i}\right)$ and
 \beq \nabla^\mathscr E e_i=f^*\left(\nabla \frac{\p}{\p y^i}\right)=\Gamma_{ij}^k(f)\frac{\p f^j}{\p x^\alpha}dx^\alpha \ts
 e_k.\eeq
We also denote by $\nabla^\mathscr E$ the extension $
\nabla^\mathscr E:\Om^p(M,\mathscr E)\>\Om^{p+1}(M,\mathscr E)$. As
a classical result, the Euler-Lagrange equation of the energy \beq
E(f)=\int_M |df|^2 dv_M\eeq is $(\nabla^\mathscr E)^* df=0$,  i.e.
\beq h^{\alpha\beta}\left(\frac{\p^2 f^i}{\p x^\alpha\p
x^\beta}-\Gamma_{\alpha\beta}^\gamma\frac{\p f^i}{\p x^\gamma}
 +\Gamma_{jk}^i\frac{\p f^j}{\p x^\alpha}\frac{\p f^k}{\p x^\beta}\right)\ts e_i=0.\label{realharmonic}\eeq
A smooth map $f:(M,h)\>(N,g)$ is called a \emph{harmonic map} if it
satisfies (\ref{realharmonic}).

\subsection{Harmonic maps from Hermitian manifolds to Riemannian manifolds}
Let $(M,h)$ be a compact Hermitian manifold, $(N,g)$ a Riemannian
manifold and $\mathscr E=f^*(TN)$ with the induced Levi-Civita
connection. Let $\{z^\alpha\}_{\alpha=1}^m$ be the local holomorphic
coordinates on $M$ where $m=\dim_\C M$.
 There are three $\mathscr E$-valued $1$-forms, namely,
\beq \bp f =\frac{\p f^i}{\p \bar z^\beta}d\bar z^\beta\ts e_i, \ \
\p f=\frac{\p f^i}{\p z^\alpha}dz^\alpha\ts e_i,\ \ \  df=\bp f+\p
f. \eeq The $\bp$-energy of a smooth map $f:(M,h)\>(N,g)$ is defined
by \beq E''(f)=\int_M |\bp f|^2 \frac{\omega_h^m}{m!} \eeq and
similarly we can define the $\p$-energy  $E'(f)$ and the total
energy $E(f)$ by \beq E'(f)=\int_M |\p f|^2 \frac{\omega_h^m}{m!},\
\ \ \ \ E(f)=\int_M |df|^2\frac{\omega_h^m}{m!}. \eeq Since $N$ is a
Riemannian manifold, it is obvious that \beq
E(f)=2E'(f)=2E''(f).\eeq Let  $f:(M,h)\>(N,g)$ be a smooth map and
$\mathscr E=f^*(TN)$. The smooth map $f$ is  {harmonic} if it is a
critical point of the Euler-Lagrange equation of the energy $E(f)$,
i.e. $(\nabla^\mathscr E)^* df=0$, which is also equivalent to
$\bp_\mathscr E^*\bp f=\p_\mathscr E^*\p f=0$. Hence, the harmonic
map equation (\ref{realharmonic}) can be written into the following
complex version \beq h^{\alpha\bar\beta}\left(\frac{\p^2 f^i}{\p
z^\alpha\p\bar z^\beta}-2\Gamma_{\alpha\bar\beta}^{\bar
\gamma}\frac{\p f^i}{\p\bar z^\gamma}+\Gamma_{jk}^i\frac{\p f^j}{\p
\bar z^\beta}\frac{\p f^k}{\p z^\alpha}\right)\ts e_i=0\label{d}\eeq
or equivalently
$$h^{\alpha\bar\beta}\left(\frac{\p^2 f^i}{\p z^\alpha\p\bar
z^\beta}-2\Gamma_{\bar\beta\alpha}^{ \gamma}\frac{\p f^i}{\p
z^\gamma}+\Gamma_{jk}^i\frac{\p f^j}{\p \bar z^\beta}\frac{\p
f^k}{\p z^\alpha}\right)\ts e_i=0$$ where \beq
\Gamma_{\alpha\bar\beta}^{\gamma}=\frac{1}{2}
h^{\gamma\bar\delta}\left(\frac{\p h_{\alpha \bar\delta}}{\p\bar
z^\beta}-\frac{\p h_{\alpha \bar\beta}}{\p\bar
z^\delta}\right)=\Gamma_{\bar\beta\alpha}^\gamma=\bar{\Gamma_{\beta\bar\alpha}^{\bar\gamma}},\
\  \Gamma_{\alpha\beta}^\gamma=\frac{1}{2}
h^{\gamma\bar\delta}\left(\frac{\p h_{\alpha \bar\delta}}{\p
z^\beta}+\frac{\p h_{\beta \bar\delta}}{\p z^\alpha}\right).\eeq

\bdefinition  $f$ is called \emph{pluri-harmonic} if it satisfies
$\p_\mathscr E\bp f=0$, i.e. \beq \left(\frac{\p^2 f^i}{\p
z^\alpha\p\bar z^\beta}+\Gamma_{jk}^i\frac{\p f^j}{\p \bar
z^\beta}\frac{\p f^k}{\p z^\alpha}\right)dz^\alpha\wedge d\bar
z^\beta\ts e_i=0.\label{pluridefinition}\eeq
 $f$ is called \emph{Hermitian harmonic} if it satisfies $\mathrm{tr}_{\omega_h}\p_\mathscr E\bp
f=0$, i.e. \beq h^{\alpha\bar\beta}\left(\frac{\p^2 f^i}{\p
z^\alpha\p\bar z^\beta}+\Gamma_{jk}^i\frac{\p f^j}{\p \bar
z^\beta}\frac{\p f^k}{\p z^\alpha}\right)\ts e_i=0.\label{pseudo2}
\eeq

 \edefinition

 \noindent  Moreover, it is easy to see
 that
\bcorollary\label{equivalent} Let $f:(M,h)\>(N,g)$ be a smooth map
from a compact K\"ahler manifold $(M,h)$ to a Riemannian manifold
$(N,g)$. Then  Hermitian harmonic maps and harmonic maps coincide.
\ecorollary

 Let $(M,g)$ be a Riemannian manifold with Levi-Civita connection
$\nabla$. The  curvature tensor is defined as $
R(X,Y)Z=\nabla_X\nabla_YZ-\nabla_Y\nabla_XZ-\nabla_{[X,Y]}Z$ for any
$X,Y,Z\in \Gamma(M,TM)$. In the local coordinates $\{x^i\}$ of $M$,
we adopt the convention $
R(X,Y,Z,W)=g(R(X,Y)Z,W)=R_{ijk\ell}X^iY^jZ^kW^\ell$. It is easy to
see that \beq R_{ijk}^\ell=\frac{\p\Gamma_{k j}^\ell}{\p
x^i}-\frac{\p\Gamma_{k i}^{\ell}}{\p x^j}+\Gamma^{p}_{k
j}\Gamma^{\ell}_{p i}-\Gamma^{p}_{k i}\Gamma^{\ell}_{p j}
\label{rcurvature}\eeq and $R_{ijk\ell}=g_{s\ell}R_{ijk}^s$.

\blemma Let $f:M\>(N,g)$ be a smooth map from a complex manifold $M$
to a Riemannian manifold $(N,g)$. Then then $(1,1)$-part of the
curvature tensor of $\mathscr E=f^*(TN)$ is \beq \bp_{\mathscr
E}\p_{\mathscr E}+\p_{\mathscr E}\bp_{\mathscr
E}=R_{ijk}^{\ell}\frac{\p f^i}{\p z^\alpha}\frac{\p f^j}{\p \bar
z^\beta}dz^\alpha\wedge d\bar z^\beta\ts e^k\ts e_\ell.\label{key}
\eeq \elemma \bproof By straightforward computations, we have
$$\p_{\mathscr E}(e_i)=\Gamma_{ji}^k\frac{\p f^j}{\p z^\alpha}dz^\alpha \cdot e_k,\ \ \ \ \ \ \ \bp_{\mathscr E}(e_i)=\Gamma_{ji}^k\frac{\p f^j}{\p \bar z^\beta}d\bar z^\beta \cdot e_k$$
Therefore, \be \bp_{\mathscr E}\p_{\mathscr
E}(e_i)=-\left[\frac{\p\Gamma_{ji}^k}{\p x^\ell}\frac{\p
f^\ell}{\p\bar z^\beta}\frac{\p f^j}{\p
z^\alpha}+\Gamma_{ji}^k\frac{\p^2 f^j}{\p z^\alpha\p\bar
z^\beta}+\Gamma_{p\ell}^k\Gamma_{ji}^\ell \frac{\p f^j}{\p
z^\alpha}\frac{\p f^p}{\p\bar z^\beta}\right]dz^\alpha\wedge d\bar
z^\beta\cdot e_k, \ee and
 \be \p_{\mathscr E}\bp_{\mathscr
E}(e_i)=\left[\frac{\p\Gamma_{ji}^k}{\p x^\ell}\frac{\p f^\ell}{\p
z^\alpha}\frac{\p f^j}{\p \bar z^\beta}+\Gamma_{ji}^k\frac{\p^2
f^j}{\p z^\alpha\p\bar z^\beta}+\Gamma_{p\ell}^k\Gamma_{ji}^\ell
\frac{\p f^j}{\p \bar z^\beta}\frac{\p f^p}{\p
z^\alpha}\right]dz^\alpha\wedge d\bar z^\beta\cdot e_k \ee

\noindent  Hence, by using the curvature formula (\ref{rcurvature}),
 we
obtain (\ref{key}). \eproof

\noindent The following constraint equation for pluri-harmonic maps
is well-known (e.g. \cite[Lemma~1.3]{OU91}), and for readers'
convenience, we include a proof here. \blemma If $f:M\>(N,h)$ is a
pluri-harmonic map, then
 \begin{eqnarray}
 R_{ikj\ell}f^i_\alpha f^j_{\bar\beta} f^k_\gamma=0
\label{D}
  \end{eqnarray}
for any $\alpha,\beta,\gamma$ and $\ell$ where $f^i_\alpha=\frac{\p
f^i}{\p z^\alpha}$ and $f^j_{\bar \beta}=\frac{\p f^j}{\p\bar
z^\beta}$. \elemma

\bproof Recall that $\mathscr E=f^*(TN)$. $f$ is pluri-harmonic,
then $\bp_\mathscr E \p f=\p_{\mathscr E}\bp f=0.$ On the other
hand, by symmetric, it is easy to see that $\p_\mathscr E \p
f=\bp_{\mathscr E}\bp f=0$.  Indeed, $$ \bp_{\mathscr E}\bp
f=\bp_{\mathscr E}\left(\frac{\p f^j}{\p\bar z^\beta}d\bar
z^\beta\ts e_i\right)=\left(\frac{\p^2 f^i}{\p\bar z^\alpha\p\bar
z^\beta}+\Gamma_{jk}^i\frac{\p f^j}{\p \bar z^\alpha}\frac{\p
f^k}{\p\bar z^\beta}\right)\left(d\bar z^\alpha\wedge d\bar z^\beta
\right)\ts e_i=0.$$ Hence, $(\p_{\mathscr E} \bp_{\mathscr
E}+\bp_{\mathscr E} \p_{\mathscr E})\p f=0$
 and by (\ref{key}), it is equivalent to
$$R_{i j k\ell} f^i_\alpha f^j_{\bar\beta} f^k_\gamma
(dz^\alpha\wedge d\bar z^\beta \wedge dz^\gamma)=0 $$ for all
$\ell$. By using symmetry and the Bianchi identity,  we obtain
$$(R_{ijk\ell}-R_{kji\ell})f^i_\alpha f^j_{\bar\beta}
f^k_\gamma=R_{ikj\ell}f^i_\alpha f^j_{\bar\beta} f^k_\gamma=0.$$
Hence, we get (\ref{D}). \eproof

\vskip 2\baselineskip

\section{Rigidity  of  pluri-harmonic maps into Riemannian manifolds}

In this section, we shall prove Theorem \ref{main1}. Let
$f:(M,h)\>(N,g)$ be a  smooth map from a compact Hermitian manifold
$M$ to a Riemannian manifold $(N,g)$. Let  $$\mathscr
F=T^{*1,0}M\otimes f^*(TN)=T^{*1,0}M\otimes\mathscr E$$ and
$\nabla^{\mathscr F}$ be the induced connection on $\mathscr F$ by
using the \textbf{Chern connection} on $M$ and Levi-Civita
connection on $N$. We have a natural decomposition $\nabla^{\mathscr
F}=\p_{\mathscr F}+\bp_{\mathscr F}$. It is obvious that $\p f \in
\Gamma(M, \mathscr F)$.

\blemma If $f:M\>(N,h)$ is a pluri-harmonic map, then $\bp_{\mathscr
F} (\p f)=0$, i.e.\beq  \bp_{\mathscr F} (\p
f)=\left[\left(\frac{\p^2 f^i}{\p z^\alpha\p\bar
z^\beta}+\Gamma_{jk}^i\frac{\p f^j}{\p \bar z^\beta}\frac{\p f^k}{\p
z^\alpha}\right) d\bar z^\beta\right] \ts d z^\alpha \ts e_i=0
\label{pluriharmonic} \eeq \elemma \bproof Note that $\bp_{\mathscr
F} (\p f)\in \Om^{0,1}(M,\mathscr F)$. Since $\p f=\frac{\p f^i}{\p
z^\alpha}dz^\alpha\ts e_i\in \Gamma(M,\mathscr F)$ and
$\nabla^{\mathscr F}$ is induced by the Chern connection on $M$ and
Levi-Civita connection on $N$, \beq \bp_{\mathscr F}\p
f=\bp_\mathscr{F}\left(\frac{\p f^i}{\p z^\alpha}dz^\alpha\ts
e_i\right)=\left(\frac{\p^2 f^i}{\p z^\alpha\p\bar
z^\beta}+\Gamma_{jk}^i\frac{\p f^j}{\p \bar z^\beta}\frac{\p f^k}{\p
z^\alpha}\right) d\bar z^\beta \ts d z^\alpha \ts e_i. \eeq By the
definition equation (\ref{pluridefinition}) of pluri-harmonic maps,
we get $\bp_{\mathscr F}\p f=0$. \eproof

\noindent The following result is a generalization of the classical
Chern-Lu identity (\cite{Che68,Lu68}, see also
\cite[Lemma~5.1]{Yang18a}).

\bproposition Let $f:(M,h)\>(N,g)$ be a pluri-harmonic map from a
 Hermitian manifold $M$ to a Riemannian manifold $(N,g)$.
Then we have \beq \p\bp u =\la\p_{\mathscr F} \p f, \p_{\mathscr
F}\p f \ra + \left( R_{\alpha\bar \beta \gamma\bar \delta} h^{\gamma
\bar \nu} h^{\mu\bar \delta} g_{ij} f^{i }_\mu {f^{j}_{\bar\nu}} -
R_{i\ell k j}f^i_{\al}f^j_{\bbe}
\left(h^{\gamma\bar\delta}f^k_{\ga}f^\ell_{\bar\delta}\right)\right)dz^\alpha\wedge
d\bar z^\beta.\label{hessian} \eeq and \beq \Delta_h u=
|\p_{\mathscr F} \p f|^2 + \left(h^{\alpha\bar\beta}R_{\alpha\bar
\beta \gamma\bar \delta}\right) h^{\gamma \bar \nu} h^{\mu\bar
\delta} g_{ij} f^{i }_\mu {f^{j}_{\bar{\nu}}}  - R_{i\ell k
j}\left(h^{\alpha\bar\beta}f^i_{\al}f^j_{\bbe}\right)
\left(h^{\gamma\bar\delta}f^k_{\ga}f^\ell_{\bar\delta}\right),\label{laplacian}\eeq
 where $u=|\p f|^2 $ and $\Delta_h u=h^{\alpha\bar\beta}\frac{\p^2 u}{\p z^\alpha\p\bar
 z^\beta}$.
\eproposition

 \bproof By using Bochner technique, we have \be \p \bp u&=&\p\la
\bp_{\mathscr F} \p f, \p f\ra+\p\la \p f, \p_{\mathscr F}\p f
\ra\\&=&\p\la \p f, \p_{\mathscr F}\p f
\ra\\
&=& \la\p_{\mathscr F} \p f, \p_{\mathscr F}\p f \ra+\la \p f,
\bp_{\mathscr F}\p_{\mathscr F}\p f \ra\\
&=&\la\p_{\mathscr F} \p f, \p_{\mathscr F}\p f \ra+\la \p f,
\left(\p_{\mathscr F}\bp_{\mathscr F}+\bp_{\mathscr F}\p_{\mathscr
F}\right)\p f \ra
 \ee
Since $\mathscr F=T^{*1,0}M\otimes\mathscr E$, we have the curvature
formula \beq \p_{\mathscr F}\bp_{\mathscr F}+\bp_{\mathscr
F}\p_{\mathscr F}=R^{T^{*1,0}M}\ts \mathrm{Id_{\mathscr
E}}+\mathrm{Id_{T^{*1,0}M}}\ts \left(\p_{\mathscr E}\bp_{\mathscr
E}+\bp_{\mathscr E}\p_{\mathscr E}\right). \label{key2}\eeq Let's
compute the curvature term on the target manifold. By using the
expression
$$\left[\mathrm{Id_{T^{*1,0}M}}\ts \left(\p_{\mathscr E}\bp_{\mathscr
E}+\bp_{\mathscr E}\p_{\mathscr E}\right)\right]\p
f=\left(R_{ijk}^\ell f^i_\alpha f^j_{\bar\beta} f^k_\delta
dz^\alpha\wedge d\bar z^\beta\right) \cdot dz^\delta \ts e_\ell
$$
we obtain \beq \left\la \p f, \left[\mathrm{Id_{T^{*1,0}M}}\ts
\left(\p_{\mathscr E}\bp_{\mathscr E}+\bp_{\mathscr E}\p_{\mathscr
E}\right)\right]\p f \right\ra =-R_{ijk\ell} f^i_\alpha f^j_{\bar
\beta} \left(h^{\gamma\bar\delta} f^k_\gamma f^\ell_{\bar
\delta}\right) dz^\alpha\wedge d\bar z^\beta  \label{bad}\eeq (Note
that, in general, (\ref{bad}) has no positivity.)  By the constraint
equation (\ref{D}) for pluri-harmonic maps,  one has $-R_{ik\ell j }
f^i_\alpha f^j_{\bar \beta} f^k_\gamma=0$ for any
$\alpha,\beta,\gamma,\ell$, and so
$$-R_{ik\ell j } f^i_\alpha f^j_{\bar
\beta} \left(h^{\gamma\bar\delta} f^k_\gamma f^\ell_{\bar
\delta}\right)=0.$$ Then by the Bianichi identity, we get \be
\left\la \p f, \left[\mathrm{Id_{T^{*1,0}M}}\ts \left(\p_{\mathscr
E}\bp_{\mathscr E}+\bp_{\mathscr E}\p_{\mathscr E}\right)\right]\p f
\right\ra &=&-R_{ijk\ell} f^i_\alpha f^j_{\bar \beta}
\left(h^{\gamma\bar\delta} f^k_\gamma f^\ell_{\bar \delta}\right)
dz^\alpha\wedge d\bar z^\beta\\&=&- R_{i\ell k j}f^i_{\al}f^j_{\bbe}
\left(h^{\gamma\bar\delta}f^k_{\ga}f^\ell_{\bar\delta}\right)dz^\alpha\wedge
d\bar z^\beta.\ee On the other hand, it is easy to see that
$$\left\la \p f, \left[R^{T^{*1,0}M}\ts \mathrm{Id_{\mathscr
E}}\right]\p f \right\ra = R_{\alpha\bar \beta \gamma\bar \delta}
h^{\gamma \bar \nu} h^{\mu\bar \delta} g_{ij} f^{i }_\mu
{f^{j}_{\bar\nu}} dz^\alpha\wedge d\bar z^\beta.$$
  Hence,
we get  formula (\ref{hessian}) and by taking trace, one has
(\ref{laplacian}).
 \eproof

Let $(N,g)$ be a  Riemannian manifold and  $T_\C N:=TN\ts \C$ be the
complexification of the real vector bundle $TN$. We can extend the
metric $g$ and $\nabla$ to $T_\C N$ in the $\C$-linear way and still
denote them by $g$ and $\nabla$ respectively.

\bdefinition    $(N,g)$ is said to have \emph{non-positive complex
sectional curvature} if \beq R(Z,\bar W, W,\bar Z)\leq 0 \eeq for
any $Z,W\in T_\C N$.\edefinition

 \bdefinition
A Hermitian manifold $(M,h)$ is called \emph{uniformly RC-positive}
if the Chern curvature tensor of the Hermitian tangent bundle
$(T^{1,0}M,h)$ is uniformly RC-positive, i.e., at every point $p\in
M$, there exists some vector $u\in T^{1,0}_pM$ such that for any
nonzero vector $v\in T^{1,0}_pM$, one has  \beq R^{h} (u,\bar
u,v,\bar v)=R_{\alpha\bar\beta\gamma\bar\delta}u^\alpha\bar u^\beta
v^\gamma\bar v^\delta>0.\eeq
 \edefinition

 \btheorem Let
$f:M\>(N,g)$ be a pluri-harmonic map from a  compact complex
manifold to a Riemannian manifold $(N,g)$ with non-positive complex
sectional curvature. If $T^{1,0}M$ is uniformly RC-positive, then
$f$ is a constant map. \etheorem

\bproof Suppose--to the contrary--that $f$ is not a constant map,
then $u=|\p f|^2$ is not constant. Suppose $u$ attains a maximum
value at some point $p\in M$. Let $h$ be a Hermitian metric on $M$
such that
 has uniformly RC-positive curvature.  At point $p$,
there exists a nonzero vector $a=(a_\alpha)$ such that the curvature
of $(T^{1,0}M,h)$ has the property
$$R_{\alpha\bar\beta\gamma\bar\delta} a_\alpha\bar a_\beta v_\gamma \bar v_\delta>0 $$
for all $v\neq 0$. In particular, we have \beq
\sum_{i,\alpha,\beta,\gamma,\delta}R_{\alpha\bar\beta\gamma\bar\delta}
a_\alpha\bar a_\beta
f^i_\delta\bar{f^i_\gamma}>0.\label{positive}\eeq On the other hand,
let $u^i=\sum_\alpha f^i_\alpha a_\alpha$, then
$$\sum_{\alpha,\beta}
R_{i\ell k j}f^i_{\al}f^j_{\bbe}
\left(h^{\gamma\bar\delta}f^k_{\ga}f^\ell_{\bar\delta}\right)a_\alpha\bar
a_\beta =\sum_{\gamma}R_{i\ell k j}u^i \bar u^j
f^k_{\ga}f^\ell_{\bar\gamma}.$$ Since $(N,h)$ has non-positive
complex sectional curvature, we deduce \beq \sum_{\gamma}R_{i\ell k
j}u^i \bar u^j f^k_{\ga}f^\ell_{\bar\gamma}\leq 0
\label{nonnegative}\eeq By formulas (\ref{hessian}),
 (\ref{positive}) and (\ref{nonnegative}), we get a contradiction.
 Hence $f$ must be a constant.
 \eproof

\noindent By using similar arguments for the equation
(\ref{laplacian}), we can also prove

\bcorollary Let $f:M\>(N,g)$ be a pluri-harmonic map from a  compact
complex manifold to a Riemannian manifold $(N,g)$ with non-positive
complex sectional curvature. If $M$ has a Hermitian metric $\omega$
such that the second Chern-Ricci curvature
$\left(h^{\alpha\bar\beta}R_{\alpha\bar \beta \gamma\bar
\delta}\right)$ is positive definite, then $f$ is a constant map.
\ecorollary

\vskip 2\baselineskip

\section{Rigidity of Hermitian harmonic maps into Riemannian
manifolds}

In this section, we prove Theorem \ref{main2} by using Theorem
\ref{main1} and  ideas in \cite{Siu80, Sam85, Sam86, JY91, JY93,
LY14}.

 \blemma\label{realidentity} Let $f:(M,h)\>(N,g)$ be a
smooth map from a compact Hermitian manifold $(M,h)$ to a Riemannian
manifold $(N,g)$. Then\beq \p\bp\{\bp f, \bp f\}
\frac{\omega_h^{m-2}}{(m-2)!}= 4\left(|\p_\mathscr E\bp
f|^2-|\mathrm{tr}_{\omega_h} \p_\mathscr E\bp
f|^2\right)\frac{\omega^m_h}{m!}+Q\cdot\frac{\omega_h^m}{m!},
\label{siu15}\eeq where $$Q=-4 h^{\alpha\bar\delta}
h^{\gamma\bar\beta}R_{ijk\ell} \frac{\p f^i}{\p z^\alpha}\frac{\p
f^k}{\p \bar z^\beta}\frac{\p f^j}{\p z^\gamma}\frac{\p f^\ell}{\p
\bar z^\delta}.$$

  \elemma \bproof It is  well-known (e.g.
\cite[Lemma~6.10]{LY14} or \cite{OU91})  by using a standard Bochner
type calculations. For readers' convenience, we sketch a proof here
following Siu's $\p\bp$ trick. It is easy to see that $\bp_{\mathscr
E}\bp f=0$. Hence, we have \beq \p\bp\{\bp f, \bp f\}=-\{\p_\mathscr
E\bp f,\p_\mathscr E\bp f\}+\{\bp f, \left( \bp_{\mathscr
E}\p_{\mathscr E}+\p_{\mathscr E}\bp_{\mathscr E}\right) \bp
f\}.\eeq Moreover, \beq-\{\p_\mathscr E\bp f,\p_\mathscr E\bp
f\}\frac{\omega_h^{m-2}}{(m-2)!}=4\left( |\p_\mathscr E\bp
f|^2-|\mathrm{tr_{\omega_h}}\p_\mathscr E\bp
f|^2\right)\frac{\omega_h^m}{m!}, \label{siu2}\eeq

\noindent and by \cite[Corollary~1.2.28]{Huy05} \beq \{\bp f, \left(
\bp_{\mathscr E}\p_{\mathscr E}+\p_{\mathscr E}\bp_{\mathscr
E}\right) \bp f\} \frac{\omega_h^{m-2}}{(m-2)!}=4\left\la \sq \left[
\bp_{\mathscr E}\p_{\mathscr E}+\p_{\mathscr E}\bp_{\mathscr
E},\Lambda_{\omega_h}\right]\bp f, \bp
f\right\ra\frac{\omega_h^m}{m!}. \eeq Moreover, \beq  \left\la \sq
\left[ \bp_{\mathscr E}\p_{\mathscr E}+\p_{\mathscr E}\bp_{\mathscr
E},\Lambda_{\omega_h}\right]\bp f, \bp
f\right\ra=-h^{\alpha\bar\delta} h^{\gamma\bar\beta}R_{ijk\ell}
\frac{\p f^i}{\p z^\alpha}\frac{\p f^k}{\p \bar z^\beta}\frac{\p
f^j}{\p z^\gamma}\frac{\p f^\ell}{\p \bar z^\delta}.  \eeq
Therefore, the formula (\ref{siu15}) follows.
 \eproof

\noindent The following formulation is essentially obtained in
\cite[Theorem~6.11]{LY14}.

\btheorem\label{pluri} Let $(M,h)$ be a compact  astheno-K\"ahler
manifold (i.e. $\p\bp\omega^{m-2}_h=0$) and $(N,g)$ a  Riemannian
manifold. Let $f:(M,h)\>(N,g)$ be a Hermitian harmonic map. If
$(N,g)$ has
 non-positive complex sectional  curvature, then $f$ is pluri-harmonic. \bproof If $f$
is Hermitian harmonic, i.e., $\mathrm{tr}_{\omega_h}\p_\mathscr E\bp
f=0$, by formula (\ref{siu15}),\beq \int_M\p\bp\{\bp f, \bp f\}
\frac{\omega_h^{m-2}}{(m-2)!}=4\int_M |\p_\mathscr E\bp
f|^2\frac{\omega^m_h}{m!}+\int_MQ\cdot\frac{\omega_h^m}{m!}. \eeq
From integration by parts, we obtain
 $$4\int_M |\p_\mathscr E\bp
f|^2\frac{\omega^m_h}{m!}+\int_MQ\cdot\frac{\omega_h^m}{m!} =0.$$ If
$(N,g)$ has  non-positive complex sectional curvature, then $Q\geq
0$ and so $\p_\mathscr E\bp f=0$, i.e. $f$ is pluri-harmonic.
 \eproof \etheorem

\noindent In particular, one has the following result obtained in
\cite[Theorem~1]{Sam86} (see also \cite[Lemma~4.2]{JY91} and
\cite[Theorem]{YZ91})

\bcorollary\label{pluri2} Let $(M,h)$ be a compact K\"ahler manifold
and $(N,g)$ a Riemannian manifold. Let $f:(M,h)\>(N,g)$ be a
harmonic map. If $(N,g)$ has non-positive complex sectional
curvature, then $f$ is pluri-harmonic. \ecorollary

\bproof It follows by Corollary \ref{equivalent} and Theorem
\ref{pluri}. \eproof

\vskip 1\baselineskip

\noindent\emph{The proof of Theorem \ref{main2}.} Let $(M,h)$ be a
compact astheno-K\"ahler manifold and $(N,g)$ be a Riemannian
manifold with non-positive complex sectional curvature. By Theorem
\ref{pluri}, every Hermitian harmonic map $f:(M,h)\>(N,g)$ is
pluri-harmonic. Let $h_0$ be a Hermitian metric on $M$ such that
$(T^{1,0}M,h_0)$ is uniformly RC-positive, then by Theorem
\ref{main1}, the pluri-harmonic map $f:M\>(N,g)$ is constant. \qed

\end{document}